 \documentclass[a4paper,12pt]{article}

 \usepackage{float}
 \usepackage{color}
 \usepackage[pdftex]{graphicx}
 \usepackage{epsfig}       
 \usepackage{subfigure}    
 \usepackage{latexsym}
 \usepackage{graphicx}
 \usepackage{epstopdf}
 \usepackage{amsfonts}
 \usepackage{mathrsfs}
 \usepackage{amsfonts}
 \usepackage{amssymb,amsthm,amsmath,amscd,fancybox,graphicx}
 \usepackage{amssymb,amsmath,euscript}
 \usepackage{amsmath}
 \usepackage{nccmath}
 \usepackage{color}
 \usepackage{fancyhdr}
 \usepackage{multirow}
 \usepackage[figuresright]{rotating}
 \usepackage{empheq}
 \usepackage{cases}
 \usepackage{indentfirst}
 \usepackage[numbers,sort&compress]{natbib}
 
 \usepackage{bm}

 \def \beginproof{\par\noindent {\bf Proof.}\ \ }
 \def \endproof{\hskip .5cm $\Box$ \vskip .5cm}

 \newcommand{\M}{\mathcal{M}}
 \newcommand{\I}{\mathcal{I}}
 \newcommand{\D}{\mathcal{D}}
 \newcommand{\X}{\mathcal{X}}

 \newcommand{\R}{\mathbb{R}}

 \DeclareMathOperator{\intt}{int}
 \DeclareMathOperator{\obj}{obj}
 
 \newtheorem{proposition}{Proposition}[section]
 
 \newtheorem{lemma}{Lemma}[section]
 \newtheorem{theorem}{Theorem}[section]
 \newtheorem{corollary}{Corollary}[section]
 
 \theoremstyle{definition}
 \newtheorem{definition}{Definition}[section]
 \newtheorem{example}{Example}[section]
 
 \newtheorem{remark}{Remark}[section]

\begin{document}

 \medskip

 \begin{center}
 {\large \bf Interval optimization problems on Hadamard manifolds: Solvability and Duality}

 \vskip 0.6cm

 Le Tram Nguyen
 \footnote{E-mail: letram07st@gmail.com.} \\
 Faculty of Mathemactics,\ \\
 The University of Da Nang
 - University of Science and Education \\
 and \\
 Department of Mathematics  \\
 National Taiwan Normal University \\
 Taipei 11677, Taiwan

 \vskip 0.6cm

 Yu-Lin Chang \footnote{E-mail: ylchang@math.ntnu.edu.tw.} \\
 Department of Mathematics  \\
 National Taiwan Normal University \\
 Taipei 11677, Taiwan

 \vskip 0.6cm

 Chu-Chin Hu \footnote{E-mail: cchu@ntnu.edu.tw.} \\
 Department of Mathematics  \\
 National Taiwan Normal University \\
 Taipei 11677, Taiwan

 \vskip 0.6cm

 Jein-Shan Chen
 \footnote{Corresponding author. E-mail: jschen@math.ntnu.edu.tw. The research is
 supported by Ministry of Science and Technology, Taiwan.} \\
 Department of Mathematics  \\
 National Taiwan Normal University \\
 Taipei 11677, Taiwan

 \vskip 0.6cm

 October 12, 2022.
\end{center}
 \noindent
 {\bf Abstract.} In this paper, we will study about the solvability and duality for interval optimization problems on Hadamard manifolds. It include the KKT conditions,  Wofle dual problem with weak duality and strong duality. These results are the complement for the solvability of interval optimization problems on Hadamard manifolds.\

 \medskip


 \medskip

 \section{Introduction}
Optimization problems have a lot applications in many research fields. Due to the image of objective functions, we have different types of problems: desterministic problems, stochastic problems, or interval problems. The last type of problem mean the value of objective functions are closed interval in $\R$. In addition, because the variation bound of the uncertain variables can be obtained only through small amount of uncertainly information, the interval programming can easily handle some optimization problems. See \cite{KS00, DK94} for more information.\par

 \medskip

A Riemannian manifold has no linear structure, in fact, it is locally  identified with Euclidean space. In this setting, the Euclidean metric is replaced by Riemannian metric and the line segments are replaced by minimal geodesics. Then, the Riemannian optimization problem is, at least locally, equivalent to the smoothly constrained optimization problem on Euclidean space. And solving nonconvex constrained problem on $\R^n$ may be equivalent to solving the convex uncostrained problem on Riemannian manifold.\par

 \medskip

There are so many optimization problems, which can not be solved on Euclidean space and require the Hadamard manifolds structure. For example, in engineering (see \cite{B14}), in controlles themornuclear fusion research (see \cite{DK94}). Recently, there are some types and methods of optimization problem have been extended from $\R^n$  to Riemannian manifolds, particularly to Hadamard manifolds \cite{AMS08, CH15, BFO10, BM12, C22, H13, B20}.\par

 \medskip
Consider the primal problem as
\begin{align}
&\min f(x)\nonumber\\
&\text{ s.t. }x\in X,  g_i(x)\le 0, i=1, ..., r\nonumber
\end{align}
where $f, g_i: \R^n\longrightarrow \R, i=1, ..., r$ are given functions.\\
The Lagrange function is defined as $L(x, \mu)=f(x)+\sum\limits_{i=1}^{r}\mu_ig_i(x)$.\\
And the dual function  $q(\mu)=\inf\limits_{x\in X} L(x, \mu)$, we have the \textit{Lagrange dual problem}
\begin{align}
&\max q(\mu)\nonumber\\
&\text{ s.t. } \mu\ge 0\nonumber
\end{align}
However, in the case of interval valued function, $q(\cdot)$ does not exist.\\
Suppose $f, g_i, i=1, ..., r$ are convex, diffirentiable functions. We have the \textit{Wolfe dual problem}
\begin{align}
&\max L(x, \mu)\nonumber\\
&\text{ s.t. } \nabla L(x, \mu)=0\nonumber
\end{align}
 \medskip

Wu (see \cite{W08}) presented about Wofle duality for interval valued optimization problems (IOPs) based on $H$-difference and $H$-derivative. However, $H$-difference  does not exist for all pair of intervals. It means the class of $H$-differentiable interval valued function is very small. In this paper, besides the extendtion of IOPs to Riemannian interval valued optimization problems (RIOPs), we also use the general Hukuhara difference ($gH$-difference) and $gH$-diffirentibility of interval valued functions on Hadamard manifolds to study about RIOPs and its dual problems. The KKT conditions of RIOPs will be studied. It is more general then the concept in \cite{C22} since the constraint is defined by the interval valued functions.

 \medskip

 \section{Premilinaries}
In this section, we will recall some knowlege about interval, computing and the oder on the set of all intervals. The basic information about Riemannian manifolds, particularly the Hadamard manifolds, is necessary for the convinience of reader. Riemannian interval value functions (RIVF) and their properties are also  important for the next sections.\par

\medskip

 Following the notations used in \cite{DK94}, let $\mathcal{I}(\mathbb{R})$ be the set
 of all closed, bounded interval in $\mathbb{R}$, i.e.,
 \[
 \mathcal{I}(\mathbb{R})=\{[\underline{a}, \overline{a}] \, | \, \underline{a},
 \overline{a}\in\mathbb{R}, \, \underline{a} \leq \overline{a}\}.
 \]
 The Hausdorff metric $d_H$ on $\mathcal{I}(\mathbb{R})$ is defined by
 \[
 d_H(A,B)=\max\{|\underline{a}-\underline{b}|, |\overline{a}-\overline{b}| \}, \quad
 \forall A=[\underline{a}, \overline{a}], \
 B=[\underline{b}, \overline{b}]\in\I(\mathbb{R}).
 \]
 Then, $(\mathcal{I}(\mathbb{R}), d_H)$ is a complete metric space, see \cite{LBD05}.
 The Minkowski sum and scalar multiplications is given respectively by
 \begin{eqnarray*}
 A+B &=& [\underline{a}+\underline{b}, \overline{a}+\overline{b}], \\
 \lambda A &=&
 \begin{cases}
 [\lambda\underline{a}, \lambda\overline{a}]& \text{ if } \lambda \geq 0, \\
 [\lambda\overline{a}, \lambda\underline{a}]& \text{ if } \lambda< 0.
 \end{cases}
 \end{eqnarray*}
 where $A=[\underline{a}, \overline{a}]$, $B=[\underline{b}, \overline{b}]$.
 Note that, $A-A=A+(-1)A \neq 0$.
 A crucial concept in achieving a useful working definition of derivative for
 interval valued functions is trying to derive a suitable difference between
 two intervals.

 \medskip

 \begin{definition}[$gH$-difference of intervals \cite{S08}]
 Let $A, B\in\mathcal{I}(\mathbb{R})$. The $gH$-difference between $A$ and $B$
 is defined as the interval $C$ such that
 \[
 C=A-_{gH}B \quad \Longleftrightarrow \quad
 \begin{cases}
 A =B+C \\
 \text{or} \\
 B =A-C.
 \end{cases}
 \]
 \end{definition}

 \medskip

 \begin{proposition}{\cite{S08}}
 For any two intervals $A=[\underline{a}, \overline{a}]$, $B=[\underline{b}, \overline{b}]$,
 the $gH$-difference $C=A-_{gH}B$ always exists and
 \[
 C =\left[ \min\{\underline{a}-\underline{b},\overline{a}-\overline{b} \}, \
 \max\{\underline{a}-\underline{b},\overline{a}-\overline{b} \} \right].
 \]
 \end{proposition}

 \medskip

 \begin{proposition}{\cite{LBD05}} \label{property-dH}
 Suppose that $A, B, C\in\I(\mathbb{R})$. Then, the following properties
 hold.
 \begin{description}
 \item[(a)] $d_H(A, B)=0$ if and only if $A=B$.
 \item[(b)] $d_H(\lambda A, \lambda B)=|\lambda|d_H(A, B)$, for all $\lambda\in\mathbb{R}$.
 \item[(c)] $d_H(A+C, B+C)=d_H(A, B)$.
 \item[(d)] $d_H(A+B, C+D) \leq d_H(A, C)+d_H(B, D)$.
 \item[(e)] $d_H(A, B)= d_H(A-_{gH}B, 0)$.
 \item[(f)] $d_H(A-_{gH}B, A-_{gH}C)=d_H(B-_{gH}A, C-_{gH}A)=d_H(B, C)$.
 \end{description}
 \end{proposition}

 Notice that, for all $A\in \I(\mathbb{R})$, we define $||A||:=d_H(A, 0)$, then $||A||$
 is a norm on $\I(\mathbb{R})$ and $d_H(A, B)=||A-_{gH}B||$.\par 
 \medskip
 There is no natural ordering on $\mathcal{I}(\mathbb{R})$, therefore we need to define it.

 \begin{definition}\label{ordering} \cite{W08}
 Let $A=[\underline{a}, \overline{a}]$ and $B=[\underline{b}, \overline{b}]$ be two
 elements of $\mathcal{I}(\mathbb{R})$. We write $A\preceq B$ if
 $\underline{a}\le \underline{b}$ and $\overline{a}\le \overline{b}$. We write
 $A\prec B$ if  $A\preceq B$ and $A\ne B$.
 Equivalently, $A\prec B$ if and only if one of the following cases holds:
 \begin{itemize}
 \item $\underline{a}<\underline{b}$ and $\overline{a}\le \overline{b}$.
 \item $\underline{a}\le\underline{b}$ and $\overline{a}< \overline{b}$.
 \item $\underline{a}<\underline{b}$ and $\overline{a}< \overline{b}$.
 \end{itemize}
 We write, $A\nprec B$ if none of the above three cases hold. If  neither
 $A\prec B$ nor $B\prec A$, we say that none of $A$ and $B$ dominates the other.\\
Let $\mathcal{A}$ and $\mathcal{B}$ be two sets of closed intervals. We write $\mathcal{A}\preceq \mathcal{B}$ if and only if $A\preceq B$ for any $A\in\mathcal{A}$ and $B\in\mathcal{B}$.
 \end{definition}

 \medskip

 \begin{lemma} \label{property-sets}
 For any elements $A, B, C$ and $D$ of $\mathcal{I}(\mathbb{R})$, there hold
 \begin{description}
 \item[(a)] $A\preceq B \ \Longleftrightarrow \ A -_{gH}B\preceq \textbf{0}$.
 \item[(b)] $A\nprec B \ \Longleftrightarrow \ A-_{gH}B\nprec \textbf{0}$.
 \item[(c)] $A\preceq B \ \Longrightarrow \ A-_{gH}C\preceq B-_{gH}C$.
 \item[(d)] If $0\preceq A$ then $A-_{gH}B\preceq C \ \Longrightarrow -B\preceq C$.
 \item[(e)] $C\preceq A-_{gH}B \ \Longleftrightarrow -A -_{gH}(-B)\preceq (-C)$.
 \end{description}
 \end{lemma}
\beginproof
The proofs of (a), (b), (c) can be found in \cite{NCHC21}\par
(d) We have $-B=[-\overline{b}, -\underline{b}]$, and
\[
-\underline{b}\le \underline{a}-\underline{b}\Longrightarrow -\underline{b}\le\max\{\underline{a}-\underline{b}, \overline{a}-\overline{b}\}\le \overline{c},
\]
\[
\begin{cases}
-\overline{b}&\le\overline{a}-\overline{b}\\
-\overline{b}&\le -\underline{b}\le\underline{a}-\underline{b}
\end{cases}\Longrightarrow -\overline{b}\le\min\{\underline{a}-\underline{b}, \overline{a}-\overline{b}\}\le \underline{c}.
\]
Then, $-B\preceq C$.\par
(e), Assume $A=[\underline{a}, \overline{a}], B=[\underline{b}, \overline{b}]\in\I(\R)$, we will proof that
\[
-(A-_{gH}B)=-A-_{gH}(-B).
\]
Infact,
\begin{align}
-A-_{gH}(-B)&=[\min\{-\overline{a}-(-\overline{b}), -\underline{a}-(-\underline{b})\}, \max\{-\overline{a}-(-\overline{b}), -\underline{a}-(-\underline{b})\}]\nonumber\\
&=[\min\{-\overline{a}+\overline{b}, -\underline{a}+\underline{b}\}, \max\{-\overline{a}+\overline{b}, -\underline{a}+\underline{b}\}]\nonumber\\
&=[\min\{-(\overline{a}-\overline{b}), -(\underline{a}-\underline{b})\}, \max\{-(\overline{a}-\overline{b}), -(\underline{a}-\underline{b})\}]\nonumber\\
&=[-\max\{\underline{a}-\underline{b}, \overline{a}-\overline{b}\}, -\min\{\underline{a}-\underline{b}, \overline{a}-\overline{b}\}]\nonumber\\
&=-[\min\{\underline{a}-\underline{b}, \overline{a}-\overline{b}\}, \max\{\underline{a}-\underline{b}, \overline{a}-\overline{b}\}]\nonumber\\
&=-(A-_{gH}B).\nonumber
\end{align}
Otherwise, we can easily to see that

\[
A\preceq B\Leftrightarrow -B\preceq -A, \qquad \forall A, B\in \I(\R).
\] 
Hence, we have
\[
C\preceq A-_{gH}B \ \Longleftrightarrow -A -_{gH}(-B)\preceq (-C), \qquad \forall A, B, C\in \I(\R).
\]
\endproof
\medskip
Before study about Riemannian interval valued functions, we need some basic knowlege about Riemannian manifold, 
which can be found  in some textbooks about Riemannian geometry, such as \cite{P16, D92, J05} and the reference therein.
 Let $\M$ be a Riemannian manifold, we denote by $T_x\M$ the tangent space of $\M$
 at $x\in \M$, and the tangent bundle of $\M$ is denoted by $T\M=\cup_{x\in\M}T_x\M$.
 For every $x, y\in\M$, the Riemannian distance $d(x,y)$ on $\M$ is defined by the
 minimal length over the set of all piecewise smooth curves joining $x$ to $y$.
 Let $\nabla$ is the Levi-Civita connection on Riemannian manifold $\M$,
 $\gamma: I\subset \mathbb{R}\longrightarrow\M$ is a smooth curve on $\M$, a vector
 field $X$ is called parallel along $\gamma$ if $ \nabla_{\gamma'}X=0$, where
 $\gamma'=\dfrac{\partial \gamma(t)}{\partial t}$. We say that $\gamma$ is a geodesic
 if $\gamma'$ is parallel along itself, in this cases $\lVert \gamma'\rVert$ is 
 constant. When $\lVert \gamma'\rVert=1, \gamma$ is said to be normalized. A geodesic
 joining $x$ to $y$ in $\M$ is called minimal if its length equals $d(x,y)$.

 \medskip

 For any $x\in\M$, let $V$ be a neighborhood of $0_x\in T_x\M$, the exponential
 mapping $\exp_{x}:V \longrightarrow \M$ is defined by $\exp_x(v)=\gamma(1)$ where
 $\gamma$ is the geodesic  such that $\gamma(0)=x$ and $\gamma'(0)=v$.
 It is known that the derivative of $\exp_x$ at $0_x\in T_x\M$ is the identity map;
 furthermore, by the Inverse Theorem, it is a local diffeomorphism. The inverse map
 of $\exp_x$ is denoted by $\exp_x^{-1}$. A Riemannian manifold is complete if for
 any $x\in\M$, the exponential map $\exp_x$ is defined on $T_x\M$.
 A simply connected, complete Riemannian manifold of nonpositive sectional curvature
 is called a Hadamard manifold. If $\M$ is a Hadamard manifold, for all $x, y\in\M$,
 by the Hopf-Rinow Theorem and Cartan-Hadamard Theorem (see \cite{J05}), $\exp_x$
 is a diffeomorphism and there exists a unique normalized geodesic joining $x$ to
 $y$, which is indeed a minimal geodesic. 

\medskip

 Let $\D\subseteq \M$ be a nonempty set, a mapping
 $f: \D\longrightarrow \mathcal{I}(\mathbb{R})$ is called a Riemannian interval
 valued function (RIVF). We write $f(x)=[\underline{f}(x), \overline{f}(x)]$ where
 $\underline{f}, \overline{f}$ are real valued functions satisfy
 $\underline{f}(x)\le \overline{f}(x)$, for all $x\in\M$. Since $\mathbb{R}^n$
 is a Hadamard manifold, an interval valued function (IVF)
 $f:U\subseteq\mathbb{R}^n\longrightarrow \mathcal{I}(\mathbb{R})$ is also a RIVF. Furthermore, since $\R\subset \I(\R)$, then a Riemannian real valued function $f: \D\longrightarrow \R$ is also a RIVF.
 \medskip

 \begin{definition}\cite{W07}
 Let $U\subseteq\mathbb{R}^{n}$ be a convex set. An IVF
 $f: U\longrightarrow\mathcal{I}(\mathbb{R})$ is said to be convex on $U$ if
 \[
 f(\lambda x_1+(1-\lambda)x_2)\preceq \lambda f(x_1)+(1-\lambda)f(x_2),
 \]
 for all $x_1, x_2 \in U$ and $\lambda\in [0, 1]$.
 \end{definition}

 \medskip

 \begin{definition} \cite{GCMD20}
 Let $U\subseteq\mathbb{R}^n$ be a nonempty set. An IVF
 $f:U\longrightarrow \mathcal{I}(\mathbb{R})$ is said to be
 \textit{monotonically increasing} if for all $x, y\in U$ there has
 \[
 x \leq y \  \Longrightarrow \  f(x)\preceq f(y).
 \]
 The function $f$ is said to be \textit{monotonically decreasing} if for all
 $x, y\in U$ there has
 \[
 x \leq y \  \Longrightarrow \  f(y)\preceq f(x).
 \]
 \end{definition}

 It is clear to see that if an IVF is monotonically increasing (or monotonically
 decreasing), if and only if both the real valued functions $\underline{f}$ and $\overline{f}$
 are monotonically increasing (or monotonically decreasing).

 \medskip

 \begin{definition}\cite{GCMD20}
 Let $\D\subseteq\M$ be a nonempty set. A RIVF
 $f:\D\longrightarrow \mathcal{I}(\mathbb{R})$ is said to be \textit{bounded below}
 on $\D$ if there exists an interval $A\in\mathcal{I}(\mathbb{R})$ such that
 \[
 A\preceq f(x), \quad \forall x\in \D.
 \]
 The function $f$ is said to be \textit{bounded above} on $\D$ if there exists an
 interval $B\in\mathcal{I}(\mathbb{R})$ such that
 \[
 f(x)\preceq B, \quad \forall x\in \D.
 \]
 The function $f$ is said to be \textit{bounded} if it is both bounded below and
 above.
 \end{definition}
 It is easy to verify that if a RIVF $f$ is bounded below (or bounded above) if and only if
 both the real-valued functions $\underline{f}$ and $\overline{f}$ are bounded
 below (or bounded above).

 \medskip

 \begin{definition}\label{geodesically_convex_set}\cite{NCHC21}
 Let $\mathcal{D}\subseteq \M$ be a geodesically convex set and 
 $f:\D\longrightarrow\I(\mathbb{R})$ be a RIVF. $f$ is called geodesically 
 convex on $\D$ if
 \[
 f(\gamma(t))\preceq (1-t) f(x)+tf(y), \quad \forall x, y\in \mathcal{D}
 \ {\rm and} \ \forall t\in [0, 1],
 \]
 where $\gamma:[0, 1]\longrightarrow \M$ is the minimal geodesic joining $x$ and $y$.

 \end{definition}

 \medskip


\begin{definition}\label{concave}
A RIVF $f$ is called geodesically concave if $-f$ is geodesically convex.
\end{definition}

 \medskip

 \begin{proposition}\label{convexfunction}\cite{NCHC21}
 Let $\D$ be a geodesically convex subset of  $\M$ and $f$ be a RIVF on $\D$. 
 Then, $f$ is  geodesically convex on $\D$ if and only if $\underline{f}$ and 
 $\overline{f}$ are geodesically convex on $\D$.
 \end{proposition}

 \medskip
 \begin{proposition}\cite{NCHC21}
 The RIVF $f:\D\longrightarrow\I(\mathbb{R})$ is geodesically convex if and only if
 for all $x, y\in\D$ and $\gamma:[0, 1]\longrightarrow \M$ is the minimal geodesic 
 joining $x$ and $y$, the IVF $f\circ\gamma$ is convex on $[0, 1]$.
 \end{proposition}

 \medskip
 
 \begin{lemma}\cite{NCHC21}
 If $f$ is a geodesically convex RIVF on $\D$ and $A$ is an interval, then the sublevel set
 \[
 \D^{A}=\{x\in\D: f(x)\preceq A\},
 \]
 is a geodesically convex subset of $\D$.
 \end{lemma}

 \medskip

 \begin{definition}\cite{NCHC21}
 Let $f:\M\longrightarrow \mathcal{I}(\mathbb{R})$ be a RIVF,
 $x_0\in\M, A=[\underline{a}, \, \overline{a}]\in \mathcal{I}(\mathbb{R})$. We say
 $\lim_{x\rightarrow x_0}f(x)=A$ if for every $\epsilon>0$, there exists
 $\delta>0$ such that, for all $x\in\M$ and $d(x, x_0)<\delta$, there holds
 $d_H(f(x), A)<\epsilon$.
 \end{definition}

 \medskip
 \begin{lemma}\cite{NCHC21}
 Let $f:M\longrightarrow \mathcal{I}(\mathbb{R})$ be a RIVF,
 $A=[\underline{a}, \overline{a}]\in \mathcal{I}(\mathbb{R})$. Then,
 \[
 \lim\limits_{x\rightarrow x_0}f(x)=A
 \quad \Longleftrightarrow \quad
 \begin{cases}
 \lim\limits_{x\rightarrow x_0}\underline{f}(x)=\underline{a}, \\
 \lim\limits_{x\rightarrow x_0}\overline{f}(x)=\overline{a}.
 \end{cases}
 \]
 \end{lemma}

 \begin{definition}[$gH$-continuity] \label{gH-continuity}\cite{NCHC21}
 Let $f$ be a RIVF on a nonempty open subset $\mathcal{D}$ of $\M$, $x_0\in \D$.
 The function $f$ is said to be $gH$-continuous at $x_0$ if for all $v\in T_{x_0}\M$
 with $\exp_{x_0}v\in \D$, there has
 \[
 \lim\limits_{||v||\rightarrow 0} \left( f(\exp_{x_0}(v))-_{gH} f(x_0) \right)=\textbf{0}.
 \]
 We call $f$ is $gH$-continuous on $\D$ if $f$ is $gH$-continuous at every $x\in \D$.
 \end{definition}

 \medskip
 \begin{remark}
 We point out couple remarks regarding $gH$-continuity.
 \begin{enumerate}
 \item When $\M\equiv\mathbb{R}^{n}$, $f$ become an IVF and $\exp_{x_0}(v)=x_0+v$. 
       In other words, Definition \ref{gH-continuity} generalizes the concept of 
       the $gH$-continuity of the IVF setting, see \cite{G17}.
 \item By Lemma 3.1 and Remark 3.1, we can see that $f$ is $gH$-continuous if and 
       only if $\underline{f}$ and $\overline{f}$ are continuous.
 \end{enumerate}
 \end{remark}

 \medskip
 
 \begin{theorem}\cite{NCHC21}
 Let $\D\subseteq \M$ be a geodesically convex set with nonempty interior and 
 $f: \D\longrightarrow \I(\mathbb{R})$ be a geodesically convex RIVF. Then, $f$ 
 is $gH$-continuous on $\intt\D$.
 \end{theorem}

 \begin{definition} \cite{LMWY11}
 Let $\mathcal{D}\subseteq \M$ be a nonempty open set and consider a function
 $f:\mathcal{D}\longrightarrow\mathbb{R}$. We say that $f$ has directional
 derivative at $x\in \mathcal{D}$ in the direction $v\in T_x\M$ if the limit
 \[
 f'(x,v)=\lim\limits_{t \to 0^{+}}\dfrac{f(\exp_{x}(tv))-f(x)}{t}
 \]
 exists, where $f'(x,v)$ is called the directional derivative of $f$ at $x$
 in the direction  $v\in T_x\M$. If $f$ has directional derivative at $x$ in
 every direction $v\in T_x\M$, we say that $f$ is directional differentiable
 at $x$.
 \end{definition}

 \medskip

 \begin{definition}[$gH$-directional differentiability \cite{C22}]
 Let $f$ be a RIVF on a nonempty open subset $\D$ of $\M$. The function $f$ is
 said to have $gH$-directional derivative at $x\in \D$ in direction $v\in T_x\M$,
 if there exists a closed bounded interval $f'(x,v)$ such that the limits
 \[
 f'(x,v)=\lim\limits_{t \to 0^{+}}\dfrac{1}{t}(f(\exp_x(tv))-_{gH}f(x))
 \]
 exists, where $f'(x,v)$ is called the $gH$-directional derivative of $f$ at $x$
 in the direction of $v$. If $f$ has $gH$-directional derivative at $x$ in every
 direction $v\in T_x\M$, we say that $f$ is $gH$-directional differentiable at $x$.
 \end{definition}

 \medskip

 \begin{lemma}{\cite{C22}}\label{gHderivative}
 Let $\mathcal{D}\subseteq \M$ be a nonempty open set and consider a RIVF
 $f: \mathcal{D}\longrightarrow \mathcal{I}(\mathbb{R})$. Then, $f$ has $gH$-directional
 derivative at $x\in\D$ in the direction $v\in T_x\M$ if and only if $\underline{f}$ and
 $\overline{f}$ have directional derivative at $x$ in the direction $v$. Furthermore,
 we have
 \[
 f'(x,v)=\left[\min\{(\underline{f})'(x, v), (\overline{f})'(x,v)\}, \,
 \max\{(\underline{f})'(x, v), (\overline{f})'(x,v)\}\right],
 \]
 where $(\underline{f})'(x, v)$ and $(\overline{f})'(x,v)$ are the directional derivatives
 of $\underline{f}$ and $\overline{f}$ at $x$ in the direction $v$, respectively.
 \end{lemma}

 \medskip

\begin{remark}\label{nonlinear}Let $\mathcal{D}\subseteq \M$ be a nonempty open set and consider  RIVFs
 $f, g: \mathcal{D}\longrightarrow \mathcal{I}(\mathbb{R})$. If,  $f$ and $g$ have $gH$-directional
 derivative at $x\in\D$ in the direction $v\in T_x\M$ then
\begin{enumerate}
\item $(-f)'(x,v)=-f'(x,v)$.
\item In general, $f'(x,v)+g'(x,v)\ne (f+g)'(x,v)$.
\end{enumerate}
\end{remark}
 \medskip
 
 \begin{theorem}\label{Exitence_gH-directional_derivetive}\cite{NCHC21}
 Let $\mathcal{D}\subseteq \M$ be a nonempty open geodesically convex set. If 
 $f:\D\longrightarrow \I(\mathbb{R})$ is a geodesically convex RIVF, then at 
 any $x_0\in\D$, $gH$-directional derivative $f'(x_0, v)$ exists for every 
 direction $v\in T_{x_0}\M$.
 \end{theorem}

 \medskip
 
 \begin{theorem} \label{geodesical-convexity}\cite{NCHC21}
 Let $f:\D\longrightarrow \I(\mathbb{R})$ be a $gH$-directional differentiable RIVF.
 If $f$ is geodesically convex on $\D$, then
 \[
 f'(x, \exp_{x}^{-1}y)\preceq f(y)-_{gH}f(x), \quad \forall x, y\in \D.
 \]
 \end{theorem}

 \medskip

 \begin{corollary}\cite{NCHC21}
 Let $\D\subseteq\M$ be nonempty open geodesically convex set and suppose that
 the RIVF $f:\D\longrightarrow \I(\mathbb{R})$ is $gH$-directional differentiable
 on $\D$. If $f$ is geodesically convex on $\D$, then
 \[
 f(y)\nprec f'(x, \exp_{x}^{-1}y)+f(x), \quad \forall x, y\in \D.
 \]
 \end{corollary}

 \medskip

 \begin{corollary}\label{concave_property}
 Let $f:\D\longrightarrow \I(\mathbb{R})$ be a $gH$-directional differentiable RIVF.
 If $f$ is geodesically concave on $\D$, then
 \[
 f(y)-_{gH}f(x)\preceq f'(x, \exp_{x}^{-1}y), \quad \forall x, y\in \D.
 \]
 \end{corollary}
\beginproof

Since $f$ is geodesically concave on $\D$, then $-f$ is geodesically convex on $\D$, for all $x, y\in \D$
\[
-f'(x, \exp_{x}^{-1}y)(-f)'(x, \exp_{x}^{-1}y)\preceq -f(y)-_{gH}(-f(x))=-(f(y)-_{gH}f(x)).
\]
Hence, by Lemma \ref{property-sets}(e), we have
\[
f(y)-_{gH}f(x)\preceq f'(x, \exp_{x}^{-1}y), \quad \forall x, y\in \D.
\]
\endproof

 \medskip

 \section{Sovability}
Consider the  Riemannian interval optimization problem (RIOP)
\begin{align}\label{RIOP}
&\min f(x)\\
&\text{ s.t }x\in \D, G_i(x)\preceq \textbf{0}, i=1, ..., r\nonumber
\end{align}
where $\D\subset\M, f, G_i:\D\longrightarrow \I(\R), i=1, ..., r$ are  RIVF. We denote by
\[
\X=\{x\in \M: x\in\D, G_i(x)\preceq \textbf{0}, i=1, ..., r\}
\]
the feasible set of RIOP (\ref{RIOP}). We also denote by
\[
\obj_{P}(f, \X):=\{f(x): x\in\X\}
\]
the set of all objective value of RIOP (\ref{RIOP}).
\begin{definition}\cite{NCHC21}
Consider problem (\ref{RIOP}). A feasible point $x_0$ is said to be an efficient point  of RIOP (\ref{RIOP}) if $f(x)\nprec f(x_0)$ for all feasible point $x$. In this case, $f(x_0)$ is called efficient objective value of RIOP (\ref{RIOP}). We denote by $\min(f,\X)$ the set of all efficient objective values of RIOP (\ref{RIOP}).
\end{definition}
 \medskip

\begin{proposition}\label{solvability1}\cite{NCHC21}
Consider the RIOP (\ref{RIOP}) with $f(x)=[\underline{f}(x), \overline{f}(x)]$ and $\X$ is the feasible set. Given any $\lambda_1, \lambda_2>0$, if $x_0\in\X$ is an optimal solution of the following problem
\[
\min_{x\in\X}h(x)=\lambda_1\underline{f}(x)+\lambda_2\overline{f}(x),
\]
then, $x_0$ is an efficient point of RIOP (\ref{RIOP}).
\end{proposition}
The constraint $G_i(x)=[\underline{G}_i(x), \overline{G}_{i}(x)]\preceq \textbf{0}$ eqivalent to $\underline{G}_i(x)\le \overline{G}_i(x)\le 0$. Since the objective function $f(x)=[\underline{f}(x), \overline{f}(x)]$, we can consider two corresponding scalar problem for (\ref{RIOP}) as follows:
\begin{align}\label{LRIOP}
&\min \underline{f}(x)\\
&\text{ s.t }x\in \D, \underline{G}_i(x)\le 0, i=1, ..., r\nonumber
\end{align}
and 
\begin{align}\label{URIOP}
&\min \overline{f}(x)\\
&\text{ s.t }x\in \D, \overline{G}_i(x)\le 0, i=1, ..., r\nonumber
\end{align}

  \medskip

 \begin{proposition}\label{solvability2}\cite{NCHC21}
 Consider problem (\ref{RIOP}) and the corresponding scalar problems (\ref{LRIOP}) and (\ref{URIOP}). The followings  hold:
 \begin{description}
 \item[(a)] If $x_0\in\D$ is an optimal solution of problems (\ref{LRIOP}) and 
            (\ref{URIOP}) simultaneously , then $x_0$ is an efficient point of 
            the RIOP (\ref{RIOP}).
 \item[(b)] If $x_0\in\X$ is an unique optimal solution of problems (\ref{LRIOP}) 
            or (\ref{URIOP}) , then $x_0$ is an efficient point of the RIOP (\ref{RIOP}).
 \end{description}
 \end{proposition}
One of impotant conditions for optimization problems is Karush–Kuhn–Tucker (KKT) condition. The next part of this section, we will derive the KKT condition for RIOP (\ref{RIOP}). At first, we will give the KKT condition for the following real valued optimization on Hadamard manifolds (ROP)
\begin{align}\label{ROP}
&\min F(x)\\
&\text{ s.t }x\in \D, g_i(x)\le 0, i=1, ..., r\nonumber
\end{align}
where $F:\D\longrightarrow \R, g_i:\D\longrightarrow\R, i=1, ..., r$. Let $X=\{x\in\D: g_i(x)\le 0, i=1, ..., r\}$ be the set of feasible point to ROP (\ref{ROP}).

\medskip
\begin{theorem}{\cite{C22}}\label{KKT_ROP}
Consider the ROP (\ref{ROP}). Let $x_0\in X$. Suppose $F, g_i, i=1, ..., r$ are geodesically convex on $\D$ and directional differentiable at $x_0$. Furthermore, for every feasible point $x\in X$, there exist scalars $\mu_i\ge 0, i=1, ..., r$, such that
\[
\begin{cases}
F'(x_0, \exp_{x_0}^{-1}x)+\sum\limits_{i=1}^{m}\mu_ig'_i(x_0, \exp_{x_0}^{-1}x)\ge 0\\
\mu_ig_i(x_0)= 0, \forall i=1, ..., r.
\end{cases}
\]
Then $x_0$ is an optimal solution of ROP (\ref{ROP}).
\end{theorem}

\medskip
\begin{theorem}\label{KKT_RIOP1}
Consider the RIOP (\ref{RIOP}). Let $x_0\in\X$. Suppose $f, G_i, i=1, ..., r$ are geodesically convex on $\D$ . Furthermore, for any feasible point $x\in\X$, there exist scalars $\underline{\mu}_i, \overline{\mu}_i\ge 0, i=1, ..., r$, such that
\[
\begin{cases}
(\underline{f})'(x_0, \exp_{x_0}^{-1}x)+\sum\limits_{i=1}^{r}\underline{\mu}_i(\underline{G}_i)'(x_0, \exp_{x_0}^{-1}x)\ge 0\\
(\overline{f})'(x_0, \exp_{x_0}^{-1}x)+\sum\limits_{i=1}^{r}\overline{\mu}_i(\overline{G}_i)'(x_0, \exp_{x_0}^{-1}x)\ge 0\\
\overline{\mu}_{i}\overline{G}_i(x_0)=\underline{\mu}_i\underline{G}_{i}(x_0)=0, \forall i=1, ..., r.
\end{cases}
\]
Then, $x_0$ is an efficient point of RIPO (\ref{RIOP}).
\end{theorem}
\beginproof
Follow Proposition (\ref{convexfunction}) and Lemma (\ref{gHderivative}), the geodesically convex and $gH$-directional differentiable properties of $f, G_i, i=1, ..., r$ are eqivalent with the geodesically convex and directional differentiable properties of $\underline{f}, \overline{f}, \underline{G}_i, \overline{G}_i, i=1, ..., r$, respectively. Then, by the assumption of Theorem \ref{KKT_RIOP1}, the assumption of Theorem \ref{KKT_ROP} are hold for the URIOP (\ref{URIOP}) and LRIOP (\ref{LRIOP}). It mean $x_0$ is an optimal solution  of problems (\ref{LRIOP}) and (\ref{URIOP}) simultaneously. By Proposition \ref{solvability2}, we have $x_0$ is an efficient point of RIOP (\ref{RIOP}).
\endproof
 \begin{example}
 Let $\M=\mathbb{R}_{++}:=\{x\in \mathbb{R} \, | \, x>0\}$ be endowed with  the
 Riemannian metric given by
 \[
 \langle u, v\rangle_x=\dfrac{1}{x^2}uv, \quad \forall u, v\in T_xM\equiv \mathbb{R}.
 \]
 Then, it is known that $\M$ is a Hadamard manifold. For all $x\in\M$, $v\in T_x\M$,
 the geodesic $\gamma:\mathbb{R}\longrightarrow \M$ such that $\gamma(0)=x, \gamma'(0)=v$
 is described by
 \[
 \gamma(t)=\exp_{x}(tv)=xe^{(v/x)t} \quad {\rm and} \quad
 \exp_{x}^{-1}y=x\ln\dfrac{y}{x}, \quad \forall y\in \M.
 \]
 We consider the RIOP as below
\begin{align}
&\min f(x)\nonumber\\
&\text{ s.t. }x\in M,  G(x)\le\textbf{ 0,} \nonumber
\end{align}
where $f: \M\longrightarrow \I(\mathbb{R})$ is defined by
 \[
 f(x)=\left[x, x+\dfrac{1}{x}\right], G(x)=[\min\{0, \ln x\}, \max\{0, \ln x\}] \quad \forall x\in \M.
 \]
We have $\ln x\le 0\Longleftrightarrow x\in (0, 1]$ then the feasible set of above problem is $\X=(0, 1]$.\\
  By the Cauchy-Schwarz inequality, for all $x>0$, we have
 \begin{center}
 $x+\dfrac{1}{x}\ge 2, $ and $x+\dfrac{1}{x}=2\Leftrightarrow x=1$,
 \end{center}
 then
 \begin{center}
 $\left[x, x+\dfrac{1}{x}\right]\nprec [1, 2], \forall x>0$.
 \end{center}
We also have $G(1)=\textbf{0}$, then $x_0=1$ is an efficient point of the RIOP.\\
 \noindent
Other hand we have
\[
\underline{f}'(1, \exp_{1}^{-1}x)=\ln x, \overline{f}'(1, \exp_{1}^{-1}x)=0,
\]
\[
\underline{G}'(1, \exp_{1}^{-1}x)=\ln x, \overline{G}'(1, \exp_{1}^{-1}x)=0,
\]
for all $x\in \X$.\\
Therefore

$\begin{cases}
\underline{f}'(x_0, \exp_{x_0}^{-1}x)+\sum\limits_{i=1}^{r}\underline{\mu}_i\underline{G}'_i(x_0, \exp_{x_0}^{-1}x)< 0\\
\overline{f}'(x_0, \exp_{x_0}^{-1}x)+\sum\limits_{i=1}^{r}\overline{\mu}_i\overline{G}'_i(x_0, \exp_{x_0}^{-1}x)< 0
\end{cases}\forall x\in (0, 1).$

 \end{example}
By Example 3.1, we see that, Theorem 3.2 is only the sufficient condition.
\medskip
\begin{theorem}\label{KKT_RIOP2}
Under the same assumption of Theorem \ref{KKT_RIOP1}. If for any feasible point $x\in\X$, there exist scalar $\lambda_1, \lambda_2>0$ and $\mu_i\ge 0, i=1, ..., r$, such that
\begin{equation}\label{KKT_con2}
\begin{cases}
\lambda_1(\underline{f})'(x_0, \exp_{x_0}^{-1}x)+\lambda_2(\overline{f})'(x_0, \exp_{x_0}^{-1}x)+\sum\limits_{i=1}^{r}\mu_i\overline{G}'_i(x_0, \exp_{x_0}^{-1}x)\ge 0\\
\mu_i\overline{G}_i(x_0)=0, \forall i=1, ..., r.
\end{cases}
\end{equation}
Then, $x_0$ is an efficient point of RIOP (\ref{RIOP}).
\end{theorem}
\beginproof
Consider the Riemannian real valued problem
\begin{align}
&\min F(x)=\lambda_1\underline{f}(x)+\lambda_2\overline{f}(x)\nonumber\\
&\text{ s.t } x\in\D, \overline{G}_i(x)\leq 0, i=1, ..., r.\nonumber
\end{align}
Then, by the condition (\ref{KKT_con2}) and Theorem (\ref{KKT_ROP}), $x_0$ is an efficient point of the above problem. Since $\overline{G}_i(x)\leq 0, i=1, ..., r$, then $\underline{G}_i(x)\leq 0, i=1, ..., r$
 and $G_i(x)\preceq , i=1, ..., r$. This implies that $x_0\in \X$. By Proposition (\ref{solvability1}), we have $x_0$ is an efficient point of RIOP (\ref{RIOP}).
\endproof

\medskip

\begin{theorem}\label{KKT_RIOP3}
Under the same assumption of Theorem \ref{KKT_RIOP1}. If for any feasible point $x\in\X$, there exist scalar   $\mu_i\ge 0, i=1, ..., r$, such that
\[
\begin{cases}
\textbf{0}\preceq f'(x_0, \exp_{x_0}^{-1}x)+\sum\limits_{i=1}^{r}\mu_iG'_i(x_0, \exp_{x_0}^{-1}x)\\
\mu_iG_i(x_0)=\textbf{0}, \forall i=1, ..., r.
\end{cases}
\]
Then, $x_0$ is an efficient point of RIOP (\ref{RIOP}).
\end{theorem}
\beginproof
By the assumption of Theorem, we have
\[
\begin{cases}
(\underline{f})'(x_0, \exp_{x_0}^{-1}x)+\sum\limits_{i=1}^{r}\underline{\mu}_i(\underline{G}_i)'(x_0, \exp_{x_0}^{-1}x)&\ge \min\{(\underline{f})'(x_0, \exp_{x_0}^{-1}x), (\overline{f})'(x_0, \exp_{x_0}^{-1}x)\}\\
&+\sum\limits_{i=1}^{r}\underline{\mu}_i\min \{(\underline{G}_i)'(x_0, \exp_{x_0}^{-1}x), (\overline{G}_i)'(x_0, \exp_{x_0}^{-1}x)\}\ge 0\\
(\overline{f})'(x_0, \exp_{x_0}^{-1}x)+\sum\limits_{i=1}^{r}\overline{\mu}_i(\overline{G}_i)'(x_0, \exp_{x_0}^{-1}x)&\ge \min\{(\underline{f})'(x_0, \exp_{x_0}^{-1}x), (\overline{f})'(x_0, \exp_{x_0}^{-1}x)\}\\
&+\sum\limits_{i=1}^{r}\underline{\mu}_i\min \{(\underline{G}_i)'(x_0, \exp_{x_0}^{-1}x), (\overline{G}_i)'(x_0, \exp_{x_0}^{-1}x)\}\ge 0\\
\overline{\mu}_{i}\overline{G}_i(x_0)=\underline{\mu}_i\underline{G}_{i}(x_0)=0, \forall i=1, ..., r.
\end{cases}
\]
Hence, by Theorem (\ref{KKT_RIOP1}) we have, $x_0$ is an efficient point of RIOP (\ref{RIOP}).
\endproof

\medskip
\begin{example}
Let $\M=\R^2$ with standard metric. Then $\M$ is a flat Hadamard manifold. We consider the RIOP as below
\begin{align}
\min f(x)=&[\min\{x_1+2x_2, 2x_1+x_2\}, \max\{x_1+2x_2, 2x_1+x_2\} ]\nonumber\\
\text{ s.t. }x\in \M,  &G_1(x)=[\min\{x_1-x_2, -x_1\}, \max\{x_1-x_2, -x_1\}]\le 0, \nonumber\\
&G_2(x)=[\min\{x_2-x_1, -x_2\}, \max\{x_2-x_1, -x_2\}]\le 0, \nonumber
\end{align}
where $x=(x_1, x_2)$.\par
It is easy to see that the feasible point set is $\X=\{(x_1, x_2)\in\R^2: 0\le x_1=x_2\}$. At $x_0=(0, 0)$, for all $x=(x_1, x_2)\in\X$, we have
\[
f'(x_0, \exp_{x_0}^{-1}x)=[\min\{2x_1+x_2, x_1+2x_2\}, \max\{2x_1+x_2, x_1+2x_2\}],
\]
\[ 
G_1'(x_0, \exp_{x_0}^{-1}x)=[-x_1, 0], G_2'(x_0, \exp_{x_0}^{-1}x)=[-x_2, 0].
\]

Therefore, with $\mu_1=\mu_2=1$, we have
\[
\begin{cases}
\textbf{0}\preceq f'(x_0, \exp_{x_0}^{-1}x)+\sum\limits_{i=1}^{r}\mu_iG'_i(\exp_{x_0}^{-1}x)\\
\mu_iG_i(x_0)=0, \forall i=1, ..., r.
\end{cases}.
\]
 Hence, $x=(0, 0)$ is a feasible point of this problem.
\end{example}
\medskip

We refer to the RIOP (\ref{RIOP}) as the \textit{primal problem} and we denote by $\min(f, \X)$ the set of  optimal values of (\ref{RIOP}). We define the interval valued Langrangian function for the primal problem RIOP (\ref{RIOP}) as follows:
\[
L(x, \mu)=f(x)+\sum\limits_{i=1}^{r}\mu_iG_i(x),
\]
for all $x\in\D$ and $\mu_i\ge 0$ for all $i=1, ..., r$.

\medskip
\begin{definition}\label{LM}
A vector $\mu^*=(\mu_1^*, ..., \mu_r^*)$ is said to be a \textit{Lagrange multiplier} for the primal problem RIOP (\ref{RIOP}) if
\[
\mu_i\ge 0, i=1, ..., r,
\]
and
\[
\min(f, \X)=\min\limits_{x\in\D}L(x,  \mu^*).
\]
\end{definition}

\medskip

\begin{proposition}\label{globalmini}
Let $\mu^*$ be a Lagrange multiplier of the primal problem RIOP (\ref{RIOP}). Then, $x^*$ is an efficient of the primal problem (\ref{RIOP}) if and only if $x^*$ is feasible and
\begin{equation}\label{global_minimum}
x^*\in\arg\min\limits_{x\in\D}L(x, \mu^*), \quad \mu_i^*G_i(x^*)=\textbf{0}, \forall i=1, ..., r.
\end{equation}
\end{proposition}
\beginproof
If $x^*$ is an efficient of the primal problem (\ref{RIOP}) then $x^*$ is feasible. Since $\mu_i^*$ is a Lagrange multiplier then 
\[
f(x^*)\in \min\limits_{x\in\D}L(x, \mu^*).
\]
Hence 
\begin{align}
&L(x^*, \mu^*)\nprec f(x^*)\nonumber\\
\Rightarrow&f(x^*)+\sum\limits_{i=1}^{r}\mu_i^*G_i(x^*)\nprec f(x^*).\nonumber
\end{align}
It implies that
\[
\sum\limits_{i=1}^{r}\mu_i^*G_i(x^*)\nprec \textbf{0}.
\]
Since $G_i(x^*)\preceq \textbf{0}, i=1, ..., r$ and $\mu_i^*\ge 0, i=1, ..., r$ then $\mu_i^*G_i(x^*)= \textbf{0}, i=1, ..., r$.\par
Conversely, if $x^*$ is feasible and (\ref{global_minimum}) holds then
\[
f(x^*)=f(x^*)+\sum\limits_{i=1}^{r}\mu_i^*G_i(x^*)=L(x^*, \mu^*)\in \min\limits_{x\in\D}L(x, \mu^*)=\min(f, \X),
\]
then $x^*$ is an efficient of the primal problem (\ref{RIOP}).
\endproof
\medskip

 \section{ Wolfe duality for  interval optimization problems on Hadamard manifolds}
Consider the primal problem (\ref{RIOP}). In this section, we assume further that $f, G_i, i=1, ..., r$ are geodesically convex and $gH$-directional differentiable on $\D$. Then, we write
\[
H(x, \mu, v)=f'(x, v)+\sum\limits_{i=1}^{r}\mu_iG'_i(x, v), \forall x\in\D, \mu_i\ge 0, i=1, ..., r,  v\in T_x\D.
\]
 \medskip
Note that,  $L'(x, \mu, v)\ne H'(x, \mu, v)$ in general.\\
We consider the Wolfe dual problem (WRIOPD) of primal problem (\ref{RIOP}) as follows:
\begin{align}\label{WRIOPD}
&\max L(x, \mu)\\ 
&\text{ s.t. } H(x, \mu, v)=\textbf{0}, x\in\D, \mu\ge 0, v\in T_x\D.\nonumber 
\end{align}
Denote by 
\[
Y=\{(x, \mu)\in \M\times\R^{r}_{+}: H(x, \mu, v)=\textbf{0}, \forall x\in\D, v\in T_x\D\},
\]
the set of all feasible points of WRIOPD. We also denote by
\[
\obj_D(L, Y)=\{L(x, \mu): (x, \mu)\in Y\},
\]
the set of all objective values of Wolfe dual problem (\ref{WRIOPD}).

 \medskip

\begin{definition}
$(\overline{x}, \overline{\mu})\in Y$ is said to be an efficient point of WRIOPD (\ref{WRIOPD}) if 
\[
L(\overline{x}, \overline{\mu})\nprec L(x, \mu), \forall (x, \mu)\in Y.
\]
\end{definition}

 \medskip

\begin{proposition}\label{Wgap}
Let $\D\subseteq \M$ be an open geodesically convex set. Assume $x_0, (x_1, \mu)$ be the feasible point of RIOP (\ref{RIOP}) and WRIOPD (\ref{WRIOPD}), respectively. If $f, G_i, i=1, ..., r$ are geodesically convex on $\D$ then $L(x_1, \mu)\preceq f(x_0)$.
\end{proposition}
\beginproof
Let $A_i=[\underline{a_i}, \overline{a_i}]\in\I(\R), i=1, ..., n$, we will proof that if $\sum\limits_{i=1}^{n}A_i=\textbf{0}$ then $\underline{a_i}=\overline{a_i}, i=1, ..., n$. Indeed
\begin{equation}\label{Wgap1}
\sum\limits_{i=1}^{n}A_i=\textbf{0}\Longrightarrow 0=\sum\limits_{i=1}^{n}\underline{a_i}\leq \sum\limits_{i=1}^{n}\overline{a_i}=0.
\end{equation}
If there exists $i\in \{1, ..., n\}$ such that $\underline{a_i}<\overline{a_i}$ then, $\sum\limits_{i=1}^{n}\underline{a_i}< \sum\limits_{i=1}^{n}\overline{a_i}$. Which is contradict with (\ref{Wgap1}). Hence, $\underline{a_i}=\overline{a_i}, i=1, ..., n$.\\
Since $(x_1, \mu)$ is a feasible point of  WRIOPD (\ref{WRIOPD}) then for some $v\in T_xD$ we have
\begin{align}\label{Wgap2}
&H(x_1, \mu, v)=\textbf{0}, \nonumber\\ 
\Rightarrow& f'(x_1, v)+\sum\limits_{i=1}^{r}\mu_iG'_i(x_1, v)=\textbf{0},\nonumber\\
\Rightarrow& (\underline{f})'(x_1, v)+ \sum\limits_{i=1}^{r}\mu_i(\underline{G}_i)'(x_1, v)=(\overline{f})'(x_1, v)+ \sum\limits_{i=1}^{r}\mu_i(\overline{G}_i)(x_1, v) =0.
\end{align}
$f$ is geodesically convex RIVF on $\D$ then by Proposition \ref{convexfunction} $\underline{f}, \overline{f}$ are geodesically convex RIFs on $\D$, we have
\begin{equation}\label{Wgap3}
\underline{f}(x_0)\ge \underline{f}(x_1)+(\underline{f})'(x_1, v).
\end{equation}
Since $G_i$ is geodesically convex and $\mu_i\ge 0, i=1, ..., r$ then 
\begin{equation}\label{Wgap4}
-\sum\limits_{i=1}^{r}\mu_i(\underline{G_i})'(x_1, v)\ge \sum\limits_{i=1}^{r}\mu_i\left((\underline{G_i})'(x_1)-(\underline{G_i})'(x_0)\right)\ge\sum\limits_{i=1}^{r}\mu_i(\underline{G_i})'(x_1) .
\end{equation}
Hence, from (\ref{Wgap2})-(\ref{Wgap3}), we obtained
\[
\underline{f}(x_0)\ge \underline{f}(x_1)+\sum\limits_{i=1}^{r}\mu_i\underline{G_i}'(x_1).
\]
Similar, we also have
\[
\overline{f}(x_0)\ge \overline{f}(x_1)+\sum\limits_{i=1}^{r}\mu_i\overline{G_i}'(x_1).
\]
Hence, $L(x_1, \mu)\preceq f(x_0)$.
\endproof

 \medskip
\begin{corollary}[Solvability]\label{Solvability}
Let $\D\subseteq\M$ be a nonempty open geodesically convex set. Let $f, G_i, i=1, ..., r$ be geodesically convex and $gH$-differentiable on $\D$. Assume that $x^*$ is a feasible point of primal problem (\ref{RIOP}) and $(\overline{x}, \overline{\mu})$ is a feasible point of dual problem (\ref{WRIOPD}). If $f(x^*)=L(\overline{x}, \overline{\mu})$, then $x^*, (\overline{x}, \overline{\mu})$ are efficient points of primal problem (\ref{RIOP}) and dual problem (\ref{WRIOPD}), respectively.
\end{corollary}
\beginproof
By Proposition \ref{Wgap}, we have 
\[
f(x^*)=L(\overline{x}, \overline{\mu})\preceq f(x)\Longrightarrow f(x)\nprec f(x^*)
\]
for all feasible point $x$ of primal problem (\ref{RIOP}), which says that $x^*$ is an efficient point of (\ref{RIOP}).\\
On the other hand, for all feasible point $(x, \mu)$ of dual problem (\ref{WRIOPD}), we have
\[
L(x, \mu)\preceq f(x^*)=L(\overline{x}, \overline{\mu}),
\]
or $(\overline{x}, \overline{\mu})$ is an efficient point of (\ref{WRIOPD}).
\endproof

 \medskip
\begin{corollary}[Weak Duality ] Consider the primal problem  (\ref{RIOP}), and the Duality problem (\ref{WRIOPD}). Assume that $\D\in\M$ be an nonempty open geodesically convex set and $f, G_i, i=1, ..., r$ be geodesically convex and $gH$-differentiable on $\D$. Then, we have 
\[
 A\preceq B, \forall A\in \max(L, Y), B\in\min(f, \X).
\]
\end{corollary}
\beginproof
The result follows immediately from Proposition \ref{Wgap}.
\endproof

\medskip

Similar as \cite{W08}, we have the definition of the no duality gap.
\begin{definition}\label{Nodualgap}
\begin{enumerate}
\item We say that the primal problem (\ref{RIOP}) and the dual problem (\ref{WRIOPD}) have no duality gap in the weak sense if and only if $\min(f, \X)\cap \max(L, Y)\ne\emptyset$.
\item  We say that the primal problem (\ref{RIOP}) and the dual problem (\ref{WRIOPD}) have no duality gap in the strong sense if and only if there exist $x^*\in\D, \mu^*\ge 0$ such that $f(x^*)\in\min(f, \X), L(x^*, \mu^*)\in\max(L, Y)$, and $f(x^*)=L(x^*, \mu^*)$.
\end{enumerate}

\end{definition}

From definition \ref{Nodualgap} and Proposition \ref{Solvability} we have some results.
\medskip

\begin{theorem}[Strong Duality 1] Let $\D\in\M$ be a nonempty open geodesically convex set and $f, G_i, i=1, ..., r$ be geodesically convex and $gH$-differentiable interval valued functions on $\D$. If one of the following conditions is satisfied
\begin{itemize}
\item there exist a feasible point $x^*$ of the primal problem (\ref{RIOP}) such that $f(x^*)\in \obj_D(L, Y)$,
\item there exist a feasible point $(\overline{x}, \overline{\mu})$ of the dual problem (\ref{WRIOPD}) such that $L(\overline{x}, \overline{\mu})\in  \obj_P(f, \X)$.
\end{itemize}
Then the primal problem (\ref{RIOP})  and the dual problem (\ref{WRIOPD}) have no gap in the weak sence.
\end{theorem}

\medskip


\begin{theorem}[Strong Duality 2] Let $\D\in\M$ be a nonempty open geodesically convex set and $f, G_i, i=1, ..., r$ be geodesically convex and $gH$-differentiable interval valued functions on $\D$. Suppose that, there exist $x^*\in \D, \mu^*\ge 0 $  such that $x^*, (x^*, \mu^*)$ are feasible points of the primal problem (\ref{RIOP}) and the dual problem (\ref{WRIOPD}), respectively. Futhermore, 
\begin{center}
$\sum\limits_{i=1}^{r}\mu_i^*G_i(x^*)=\textbf{0}.$
\end{center}
Then the primal problem (\ref{RIOP})  and the dual problem (\ref{WRIOPD}) have no gap in the strong sence.
\end{theorem}

 \section{Conclusions}
In this paper, we have considered the Wofle duality for interval optimization problems on Hadamard manifolds, for which we establish the weak duality and strong duality, include the case have no gap in weak sence and thecase have no gap in strong sence.  The constraint functions are interval valued, which is more general then the previous papers. The Lagrange multiplier was considered together with some KKT conditions.\par
For more general, in the future, we may either study the theory for Riemannian interval optimization problems (RIOPs) and their duality problems on the general Riemannian manifolds.

 \end{document}